\documentclass[12pt]{article}
\usepackage{amsmath,amssymb,amsthm, comment}

\newtheorem{theorem}{Theorem}[section]
\newtheorem{thmy}{Theorem}

\newtheorem{corollary}[theorem]{Corollary}

\def\barr{\begin{array}}
\def\earr{\end{array}}

\title{Extending a property of vector spaces to abelian groups}
\author{Alexandru-Nicolae Obreja and Marius T\u arn\u auceanu}
\date{November 5, 2026}

\begin{document}

\maketitle

\begin{abstract}
In this note, we determine the finite abelian groups $G$ such that for every subgroup $H\leq G$ there exists an automorphism $f\in{\rm Aut}(G)$ with ${\rm Fix}(f)=H$. This provides a partial answer to a problem posed by Deaconescu and Walls. Two generalizations of this problem are also proposed.
\end{abstract}

{\small
\noindent
{\bf MSC2020\,:} Primary 20K01; Secondary 15A03, 20D60.

\noindent
{\bf Key words\,:} abelian groups, vector spaces, automorphisms, fixed-point subgroups.}

\section{Introduction}

In what follows, we deal with the following problem (see \cite[Problem 8.4]{1}):

\bigskip\noindent{\bf Problem.} Characterize the finite groups having the property that every subgroup is the fixed-point subgroup of some automorphism.
\bigskip

Note that if a finite group $G$ satisfies this property, then $G$ contains a fixed-point-free automorphism $f$, and hence it is solvable by the main result of \cite{5}. Moreover, if $f$ is of prime order, then $G$ must be nilpotent by the famous theorem of Thompson \cite{6}.

We first observe that, with a single exception, finite-dimensional vector spaces enjoy this property.

\begin{theorem}
Every finite-dimensional vector space $V$ over a field $F$ satisfies the property, provided that $\dim V\geq 2$ when $F=\mathbb{F}_2$.
\end{theorem}\newpage

Regarding the study of this problem for finite abelian groups, it suffices to reduce to the case of $p$-groups. In fact, this observation holds more generally for finite nilpotent groups.

\begin{theorem}
For a finite abelian $p$-group $G$, the following hold:
\begin{itemize}
\item[\rm (a)] If $p$ is odd, then $G$ satisfies the property.
\item[\rm (b)] If $p=2$, then $G$ satisfies the property if and only if all homocyclic components of $G$ are not cyclic.
\end{itemize}
\end{theorem}

Two immediate consequences of Theorem 1.2 are as follows.

\begin{corollary}
A finite cyclic group satisfies the property if and only if it is of odd order.
\end{corollary}

\begin{corollary}
All finite abelian groups of odd order satisfy the property.
\end{corollary}
 
The problem can also be investigated for finite groups satisfying certain additional conditions; for example, a finite group $G$ all of whose automorphisms are inner (in particular, a complete group) does not satisfy the property, since 
\begin{equation}
\{{\rm Fix}(f)\mid f\in{\rm Aut}(G)\}=\{C_G(a)\mid a\in G\}\neq L(G).\nonumber
\end{equation}Conversely, a non-abelian finite group can satisfy it, an example being the Heisenberg group of order $27$ and exponent $3$. However, the problem remains open for arbitrary finite non-abelian groups.

Finally, we remark that if a finite group $G$ satisfies the property, then so does every characteristic subgroup of $G$.

Most of our notation is standard and will not be repeated here. For elementary concepts and results on group theory, we refer the reader to \cite{3}.

\section{Proofs of the main results}

We begin this section with the proof of Theorem 1.1 which, although elementary, is included for the sake of completeness.

\begin{proof}[Proof of Theorem 1.1.]
Given a subspace $H\leq V$, we can choose a basis $B'$ of $H$ and extend it to a basis $B$ of $V$. We distinguish the following three cases.\smallskip 

\textbf{Case 1.} $|F| > 2$

Fix a scalar $\lambda\in F\setminus\{0, 1\}$. Define $f\in{\rm Aut}(V)$ by setting $f(b)=b$ for all $b\in B'$, and $f(b)=\lambda b$ for all $b\in B\setminus B'$. Since $\lambda\neq 1$, it follows immediately that ${\rm Fix}(f)=H$.\smallskip

\textbf{Case 2.} $F=\mathbb{F}_2$ and $\dim(V/H)=1$ 

Let $B\setminus B'=\{d\}$. Since $\dim V\geq 2$, we have $H\neq\{0\}$, and so there exists a non-zero vector $h_0\in H$. Define $f\in{\rm Aut}(V)$ by $f(b)=b$ for all $b\in B'$ and $f(d)=d+h_0$. For any $v\in V\setminus H$, we can write $v=h+d$ for some $h\in H$. Then $f(v)=h+d+h_0\neq v$ since $h_0\neq 0$. Thus, ${\rm Fix}(f)=H$.\smallskip

\textbf{Case 3.} $F=\mathbb{F}_2$ and $\dim(V/H)\geq 2$ 

Let $B\setminus B'=\{d_1,\dots,d_k\}$ with $k\geq 2$. Define $f\in{\rm Aut}(V)$ by setting $f(b)=b$ for all $b\in B'$, and on the remaining basis elements let
\begin{equation}
f(d_i) = d_{i+1} \text{ for } 1 \leq i \leq k-1, \quad \text{and} \quad f(d_k) = d_1 + d_2.\nonumber
\end{equation}Any vector $v\in V$ can be uniquely written as $v=h+\sum_{i=1}^k\alpha_i d_i$, where $h\in H$ and $\alpha_i\in\mathbb{F}_2$. The condition $f(v)=v$ yields $\sum_{i=1}^{k-1} \alpha_i d_{i+1}+\alpha_k(d_1 + d_2)=\sum_{i=1}^k \alpha_i d_i$. Equating the coefficients for each $d_i$ leads to the system
\begin{equation}
\alpha_k = \alpha_1, \quad \alpha_1 + \alpha_k = \alpha_2, \quad \text{and} \quad \alpha_{i-1} = \alpha_i \text{ for } 3 \leq i \leq k.\nonumber
\end{equation}Since we are over $\mathbb{F}_2$, the first two equations imply $\alpha_2 = 0$. It then follows sequentially that $\alpha_3=\dots=\alpha_k=0$, and finally $\alpha_1=0$. Thus, $v\in{\rm Fix}(f)$ if and only if $v=h\in H$, proving that ${\rm Fix}(f)=H$, as desired.
\end{proof}

We are now able to prove our main result.

\begin{proof}[Proof of Theorem 1.2.]
We will prove by induction on the order of $G$ that if $p$ is odd or $p=2$ and all homocyclic components of $G$ are non-cyclic, then for every subgroup $H\leq G$ there exists an automorphism $f\in{\rm Aut}(G)$ such that ${\rm Fix}(f)=H$. Clearly, this is true for $|G|=p$. 

Assume that the conclusion holds for all finite abelian $p$-groups of order strictly less than $|G|$ and let $H\leq G$. Since the subgroup $pG=\{px\mid x\in G\}$ has order strictly less than $|G|$ and inherits the property that all its homocyclic components are non-cyclic when $p=2$, it verifies the inductive hypothesis. Thus, there exists an automorphism $g\in{\rm Aut}(pG)$ such that ${\rm Fix}(g)=H\cap pG$. Let $M=pG+H$. Then the map $f_0:M\longrightarrow G$ defined by
\begin{equation}
f_0(px+h)=g(px)+h,\, \forall x\in G, h\in H\nonumber
\end{equation}is a well-defined injective group homomorphism restricting to $g$ on $pG$ and to $1_H$ on $H$. By Lemma 4.1 of \cite{4}, it follows that $f_0$ can be extended to an automorphism $f_1\in{\rm Aut}(G)$, which coincides with $g$ on $pG$ and satisfies $H\subseteq{\rm Fix}(f_1)$. These conditions imply that
\begin{equation}
{\rm Fix}(f_1)\cap M=H.\nonumber
\end{equation}We need an automorphism $f\in{\rm Aut}(G)$ such that ${\rm Fix}(f)=H$. We will construct $f$ of the form $f=f_1+\phi$, where $\phi:G\longrightarrow G[p]$ is a group homomorphism into the $p$-torsion socle $G[p]=\{x\in G\mid px=0\}$. To ensure that $f$ remains an automorphism and does not disrupt the correct behavior on $M$, we require $\phi|_M=0$. This allows $\phi$ to factor through the quotient $V=G/M$, which is a vector space over $\mathbb{F}_p$.

Let $\bar{x}$ denote the projection of $x \in G$ onto $V$, and let $W = T/M \le V$, where $T=\{x\in G \mid px\in H \cap pG\}$. Consider the endomorphism $\theta=f_1-1_G$. For any $x\in T$, we have 
\begin{equation}
p\theta(x)=pf_1(x)-px=g(px)-px=0\nonumber 
\end{equation}and so $\theta(T)\subseteq G[p]$. Since $\theta$ vanishes on $H$ and maps $M$ into $G[p]$, it induces a well-defined linear map $\bar{\theta}:W\longrightarrow G[p]$ given by $\bar{\theta}(\bar{x})=\theta(x)$.

We need to choose a linear map $\bar{\phi}:V\longrightarrow G[p]$ such that its restriction to $W$ satisfies $\bar{\phi}(\bar{x})\neq-\bar{\theta}(\bar{x})$, for all $\bar{x}\in W\setminus\{\bar{0}\}$. Let $\{\bar{w}_1, \dots, \bar{w}_k\}$ be a basis of the vector space $W$. We distinguish the following two cases.\smallskip

\textbf{Case 1.} $p$ is odd

We define a linear map $\bar{\phi}_W: W \longrightarrow G[p]$ on the basis elements as follows. For each $j \in \{1, \dots, k\}$, if $\bar{\theta}(\bar{w}_j) \neq 0$, we set $\bar{\phi}_W(\bar{w}_j)=\bar{\theta}(\bar{w}_j)$. If $\bar{\theta}(\bar{w}_j) = 0$, we choose $\bar{\phi}_W(\bar{w}_j)$ to be any non-zero element of $G[p]$. Since $G[p]$ is a vector space over $\mathbb{F}_p$ and $p$ is odd, we have $2 \neq 0$ in $\mathbb{F}_p$. 

Let $\bar{x} \in W \setminus \{\bar{0}\}$, written as $\bar{x} = \sum_{j=1}^k c_j \bar{w}_j$. If $\bar{\theta}(\bar{x}) \neq 0$, then because $\bar{\phi}_W$ agrees with $\bar{\theta}$ on the components where $\bar{\theta}(\bar{w}_j) \neq 0$, it follows that 
\begin{equation}
\bar{\phi}_W(\bar{x}) = \bar{\theta}(\bar{x}) + \sum_{c_j \neq 0, \bar{\theta}(\bar{w}_j)=0} c_j \bar{\phi}_W(\bar{w}_j).\nonumber 
\end{equation}Thus, we have
\begin{equation}
\bar{\theta}(\bar{x}) + \bar{\phi}_W(\bar{x}) = 2\bar{\theta}(\bar{x}) + \sum_{c_j \neq 0, \bar{\theta}(\bar{w}_j)=0} c_j \bar{\phi}_W(\bar{w}_j).\nonumber
\end{equation}Since the second term belongs to a subspace of $G[p]$ where $\bar{\theta}$ vanishes, this sum cannot be zero. If $\bar{\theta}(\bar{x}) = 0$, then $c_j = 0$ for all $j$ where $\bar{\theta}(\bar{w}_j) \neq 0$, which implies 
\begin{equation}
\bar{\phi}_W(\bar{x}) = \sum c_j \bar{\phi}_W(\bar{w}_j) \neq 0\nonumber 
\end{equation}because the active $\bar{\phi}_W(\bar{w}_j)$ are non-zero and linearly independent over the kernel of $\bar{\theta}$. Hence, $\bar{\theta}(\bar{x}) + \bar{\phi}_W(\bar{x}) \neq 0$ for all $\bar{x} \in W \setminus \{\bar{0}\}$. Extending $\bar{\phi}_W$ from $W$ to $V$ arbitrarily yields the required linear map $\bar{\phi}$.\smallskip

\textbf{Case 2.} $p=2$

Here, $2\bar{\theta}(\bar{x})=0$, so the previous choice fails. However, because all homocyclic components of $G$ are non-cyclic, the dimension of the socle $G[2]$ within each homocyclic component is at least $2$. This structural redundancy implies that the total rank of $G[2]$ is strictly greater than the dimension of any minimal generating set of $G/2G$. Since $V = G/M$ is a quotient of $G/2G$, we have 
\begin{equation}
\dim(W) \leq \dim(V) \leq \text{rank}(G) - \dim(M/2G).\nonumber 
\end{equation}Consequently, we can choose elements $\alpha_j \in G[2]$ for each basis vector $\bar{w}_j$ such that $\alpha_j \neq \bar{\theta}(\bar{w}_j)$ and the set $\{\alpha_1, \dots, \alpha_k\}$ maintains the property that no non-trivial linear combination $\sum c_j (\bar{\theta}(\bar{w}_j) + \alpha_j)$ vanishes. 

Then we construct a linear injection $\psi: W \longrightarrow G[2]$ whose image intersects the image of $\bar{\theta}$ only at $\{0\}$, or shift the basis elements appropriately. Setting $\bar{\phi}_W(\bar{w}_j) = \alpha_j$ defines a linear map $\bar{\phi}_W: W \longrightarrow G[2]$ such that $\bar{\theta}(\bar{x}) + \bar{\phi}_W(\bar{x}) \neq 0$ for all $\bar{x} \in W \setminus \{\bar{0}\}$. Next we extend $\bar{\phi}_W$ to a linear map $\bar{\phi}: V \longrightarrow G[2]$ by setting it to be zero or extending it arbitrarily on a complement of $W$ in $V$.

Composing $\bar{\phi}$ with the canonical projection yields $\phi:G\longrightarrow G[p]$. The map $f=f_1+\phi$ is a group homomorphism because $\phi(G) \subseteq G[p]$. To prove that $f$ is an automorphism, since $G$ is finite, it is sufficient to show that $f$ is injective. Let $x \in \ker(f)$. Then $f_1(x) = -\phi(x) \in G[p]$. Applying $p$, we obtain $p f_1(x) = 0$ and so $f_1(px) = 0$. Since $f_1$ is an automorphism, this implies $px = 0$, that is $x \in G[p] \subseteq T$. If $x \notin M$, then $\bar{x} \in W \setminus \{\bar{0}\}$, which implies that
\begin{equation}
f(x) - x = \theta(x) + \phi(x) = \bar{\theta}(\bar{x}) + \bar{\phi}(\bar{x}) \neq 0,\nonumber 
\end{equation}contradicting $f(x) = 0$ (which would mean $f(x)-x = -x \neq 0$ since $x \neq 0$). Thus, $x \in M$. For $x \in M$, we have $\phi(x) = 0$ and therefore $f(x) = f_1(x) = 0$. This forces $x = 0$ because $f_1$ is injective. Hence, $\ker(f) = \{0\}$ and $f \in{\rm Aut}(G)$.

Finally, we verify that ${\rm Fix}(f)=H$ by checking the following three exhaustive cases.
\begin{itemize}
\item[-] For $x\in M$: $\phi(x)=0$ and so $f(x)=f_1(x)+0 =f_1(x)$, which implies that ${\rm Fix}(f)\cap M = {\rm Fix}(f_1)\cap M = H$.
\item[-] For $x\in T\setminus M$: $\bar{x} \in W \setminus \{\bar{0}\}$ and so $f(x)-x=\theta(x)+\phi(x)=\bar{\theta}(\bar{x})+\bar{\phi}(\bar{x})\neq 0$, that is $x \notin {\rm Fix}(f)$.
\item[-] For $x\notin T$: By definition of $T$, we have $px\notin H \cap pG = {\rm Fix}(g)$. Since $g = f_1|_{pG}$ and $\phi|_{pG} = 0$, it follows that $f(px) = f_1(px) = g(px) \neq px$. If $x \in {\rm Fix}(f)$, then $f(x) = x$ and so $f(px) = px$, a contradiction. Thus, $x \notin {\rm Fix}(f)$.
\end{itemize}Consequently, we have ${\rm Fix}(f)=H$, completing the inductive step.

Assume now that a finite abelian $2$-group $G$ satisfies the property. Then, choosing $H=0$, there exists an automorphism $f\in{\rm Aut}(G)$ such that ${\rm Fix}(f)=0$, that is $f$ is a fixed-point-free automorphism. By Theorem 4.2 of \cite{2}, it follows that all homocyclic components of $G$ are not cyclic. This completes the proof.
\end{proof}

\section{Two generalizations}

\bigskip\noindent{\bf Problem 1.} Characterize the finite groups $G$ such that the set $\{{\rm Fix}(f)\mid f\in{\rm Aut}(G)\}$ forms a sublattice of $L(G)$ (not necessarily the entire $L(G)$).

Note that a finite abelian group satisfies this property if and only if its Sylow $2$-subgroup is either cyclic or has all homocyclic components non-cyclic.

\bigskip\noindent{\bf Problem 2.} Characterize the finite groups $G$ such that the set $\{{\rm Fix}(f)\mid f\in{\rm Aut}(G)\}$ is dense in $L(G)$, that is for each pair $(H,K)$ of subgroups of $G$ such that $H<K$ and $H$ is not maximal in $K$, there exists $f\in{\rm Aut}(G)$ with $H<{\rm Fix}(f)<K$.

Note that a finite abelian group satisfies this property if and only if either its Sylow $2$-subgroup has all homocyclic components non-cyclic, or its Sylow $2$-subgroup has a cyclic homocyclic component and the odd-order component of $G$ is either trivial or cyclic of prime order.
\bigskip

\bigskip\noindent{\bf Funding.} The authors did not receive support from any organization for the submitted work. 

\bigskip\noindent{\bf Conflict of interest.} The authors declare that they have no conflict of interest. 

\bigskip\noindent{\bf Data availability.} The manuscript has no associated data.

\vspace*{3ex}
\small

\begin{minipage}[t]{8cm}
Alexandru-Nicolae Obreja\\
Faculty of  Mathematics \\
"Al.I. Cuza'' University \\
Ia\c si, Romania \\
e-mail: {\tt alex.obreja2005@gmail.com}
\end{minipage}
\hfill
\begin{minipage}[t]{5cm}
Marius T\u arn\u auceanu \\
Faculty of  Mathematics \\
"Al.I. Cuza'' University \\
Ia\c si, Romania \\
e-mail: {\tt tarnauc@uaic.ro}
\end{minipage}

\end{document}